\documentclass[11pt,leqno,oneside]{article}

\usepackage{amsmath}
\usepackage{amsthm}
\usepackage{amstext}
\usepackage{amsopn}
\usepackage{amsbsy}

\newtheorem*{main-theorem}{Main Theorem}
\newtheorem{theorem}{Theorem}

\newtheorem{lemma}{Lemma}

\begin{document}

\title{The first Dirichlet Eigenvalue of a Compact Manifold and
the Yang Conjecture
\thanks{2000 Mathematics Subject Classification Primary 58J50, 35P15; Secondary 53C21}}

\author{Jun LING}
\date{}

\maketitle

\begin{abstract}
We give a new estimate on the lower bound of the first Dirichlet
eigenvalue of a compact Riemannian manifold with negative lower
bound of Ricci curvature and provide a solution for a conjecture
of H.~C.~Yang.
\end{abstract}

\section{Introduction}\label{sec-intro}
It has been proved by P.~Li and S.~T.~Yau\cite{liy1} that if $M$
is an n-dimensional closed Riemannian manifold with Ricci
curvature Ric$(M)$ bounded below by $(n-1)\kappa$ with constant
$\kappa<0$, then the first non-zero eigenvalue $\lambda$ of the
Laplacian of $M$ has the lower bound
\[
\lambda \geq
\frac{1}{2(n-1)d^2}{\exp\{-1-\sqrt{1+4(n-1)^2d^2|\kappa|}\}},
\]
where $d$ is the diameter of $M$. H.~C.~Yang\cite{yang} improved
the above estimate to the following
\[
\lambda \geq \frac{\pi^2}{d^2}\exp\{-C_n\sqrt{(n-1)|\kappa|d^2}\},
\]
where $C_n=\max\{\sqrt{n-1}, \sqrt{2}\}$. Yang further conjectured
that
\[
\lambda \geq \frac12(n-1)\kappa + \frac{\pi^2}{d^2}.
\]
If $M$ has a boundary, H.~C.~Yang conjectured that the above
estimate holds for the first Dirichlet eigenvalue as well. In this
paper, we give a new estimate on the lower bound of the first
Dirichlet eigenvalue of an $n$-dimensional Riemannian manifold
with negative lower bound of Ricci curvature. The result provides
a solution for the conjecture of H.~C.~Yang. Let
$\textup{dist}(\cdot, \cdot)$ be the distance on $M$. We have the
following result.
\begin{theorem} \label{main-thm}
If $(M, g)$ is an $n$-dimensional compact Riemannian manifold with
boundary. Suppose that the boundary $\partial M$ of the manifold
$M$ has nonnegative mean curvature with respect to the outward
normal and that the Ricci curvature of $M$ has lower bound
\begin{equation}\label{ricci-bound}
\textup{Ric}(M)\geq (n-1)\kappa
\end{equation}
for some constant $\kappa<0$. Then the first Dirichlet eigenvalue
$\lambda$ of the Laplacian $\Delta$ of $M$ satisfies the
inequality
\[
\lambda \geq
\frac{1}{1-(n-1)\kappa/(2\lambda)}\,\frac{\pi^2}{d^2}>0
\]
and $\lambda$ has the lower bound
\begin{equation}                    \label{main-bound}
\lambda \geq \frac12(n-1)\kappa+\frac{\pi^2}{d^2},
\end{equation}
where $d$ is the diameter of the largest interior ball in $M$,
\[
d=2r \qquad\textup{and}\qquad r=\max_{x\in M} \textup{dist}(x,
\partial M).
\]
\end{theorem}

If Ric$(M)\geq (n-1)\kappa$ with constant $\kappa>0$, it is known
that the first Dirichlet eigenvalue $\lambda$ has a lower bound as
the above. Therefore the lower bound in (\ref{main-bound}) is
universal for all three cases, constant $\kappa>0, =0$ or $<0$.

In the next section, we derive some preliminary estimates and
conditions for test functions. In the last section we construct
the needed test function and prove the main result. In the proof
of the main result, instead of using the Zhong-Yang's canonical
function or the "midrange" of the normalized eigenfunction of the
first eigenvalue, we use a function $\xi$ that the author
constructed in \cite{ling1} to build the suitable test function.

\section{Preliminary Estimates}\label{sec-pre-es}
Let $v$ be a normalized eigenfunction of the first Dirichlet
eigenvalue of the Laplacian $\Delta$ such that
\begin{equation}                                \label{v-d-con}
\sup_{M}v=1, \quad \inf_{M}v=0.
\end{equation}
The function $v$ satisfies the following
\begin{equation}                                 \label{v-d-eq}
\Delta v=-\lambda v\quad \textrm{in }M
\end{equation}
\begin{equation}                                \label{v-d-boundary}
v=0\qquad \textrm{ on }\partial M.
\end{equation}

We first use gradient estimate in \cite{li}-\cite{liy1} and
\cite{sy} to derive following estimate.
\begin{lemma}       \label{pre-es-d-lemma}
The function $v$ satisfies the following
\begin{equation}                        \label{basic3-d}
\frac{\left |\nabla v\right |^2}{b^2-v^2} \leq\lambda(1+\beta),
\end{equation}
where $\beta = -(n-1)\kappa/\lambda>0$ and $b>1$ is an arbitrary
constant.
\end{lemma}

\begin{proof}\quad
Consider the function
\begin{equation}                \label{p-of-x-def}
P(x)=|\nabla v|^2+Av^2,
\end{equation}
where $A=\lambda -(n-1)\kappa+\epsilon$ for small $\epsilon>0$.
The function $P$ must achieve its maximum at some point $x_0\in
M$. We claim that
\begin{equation}\label{gra-p-eq-0}
\nabla P(x_0)=0.
\end{equation}
If $x_0\in M\backslash \partial M$, (\ref{gra-p-eq-0}) is
obviously true. Suppose that $x_0\in\partial M$. Choose a local
orthonormal frame $\{e_1, e_2,\cdots, e_{n}\}$ of $M$ about $x_0$
so that $e_n$ is the unit outward normal vector field near $x_0\in
\partial M$ and $\{e_1, e_2,\cdots, e_{n-1}\}|_{\partial M}$ is a
local frame of $\partial M$ about $x_0$. The existence of such
local frame can be justified as the following. Let $e_n$ be the
local unit outward normal vector field of $\partial M$ about
$x_0\in \partial M$ and $\{e_1, \cdots, e_{n-1}\}$ the local
orthonormal frame of $\partial M$ about $x_0$. By parallel
translation along the geodesic $\gamma(t)=\exp_{x_0}te_n$, we may
extend $e_1$, $\cdots$, $e_{n-1}$ to local vector fields of $M$.
Then the extended frame $\{e_1, e_2,\cdots, e_{n}\}$ is what we
need. Note that $\nabla_{e_n}e_i=0$ for $i\leq n-1$. Since
$v|_{\partial M}=0$, we have $v_i(x_0)=0$ for $i\leq n-1$.

$P(x_0)$ is the maximum implies that
\begin{equation}\label{p-i-eq-0}
 P_i(x_0)= 0\qquad
\textrm{for }i\leq n-1
\end{equation}
and
\begin{equation}\label{p-n-geq-0}
P_n(x_0)\geq 0.
\end{equation}
Using (\ref{v-d-con})-(\ref{v-d-boundary}) in the following
arguments, then we have that at $x_0$,
\begin{eqnarray}      \label{gra-v-leq-0}
& &{ }\nabla_{e_n}(|\nabla
v|^2)(x_0)=\sum_{i=1}^n 2v_iv_{in}=2v_nv_{nn}\nonumber\\
& &=2v_n(\Delta^M v -\sum_{i=1}^{n-1}v_{ii}) =2v_n(-\lambda
v-\sum_{i=1}^{n-1}v_{ii})\nonumber\\
&
&=-2v_n\sum_{i=1}^{n-1}v_{ii}=-2v_n\sum_{i=1}^{n-1}(e_ie_iv-\nabla^M_{e_i}e_i
v)\nonumber\\
& &=2v_n\sum_{i=1}^{n-1}\nabla^M_{e_i}e_iv
=2v_n\sum_{i=1}^{n-1}\sum_{j=1}^{n}g(\nabla^M_{e_i}e_i
,e_j)v_j\nonumber\\
& &=2v_n^2\sum_{i=1}^{n-1}g(\nabla^M_{e_i}e_i ,e_n)
=-2v_n^2\sum_{i=1}^{n-1}g(\nabla^M_{e_i}e_n ,e_i)\nonumber\\
& &=-2v_n^2\sum_{i=1}^{n-1}h_{ii}=-2v_n^2(x_0)m(x_0)\nonumber\\
& &\leq 0 \quad\textrm{by the non-negativity of }m,
\end{eqnarray}
where $g(,)$ is the Riemann metric of $M$, $(h_{ij})$ is the
second fundamental form of $\partial M$ with respect to the
outward normal $e_n$ and $m$ is the mean curvature of $\partial M$
with respect to $e_n$.

Noticing that $v|_{\partial M}=0$, we have
\begin{equation}\label{p-n-leq-0}
P_n(x_0)=\nabla_{e_n}(|\nabla v|^2)|_{x_0}+2Av(x_0)v_n(x_0) )\leq
0.
\end{equation}
Now (\ref{p-i-eq-0}), (\ref{p-n-geq-0}) and (\ref{p-n-leq-0})
imply that $P_n(x_0)=0$.

Thus (\ref{gra-p-eq-0}) holds, no matter $x_0\not\in\partial M$ or
$x_0\in\partial M$. By (\ref{gra-p-eq-0}) and the Maximum
Principle, we have
\begin{equation}\label{max-prin}
\nabla P(x_0)=0 \qquad \textrm{and}\qquad \Delta P(x_0)\leq 0.
\end{equation}

We are going to show further that $\nabla v(x_0)=0$. If on the
contrary, $\nabla v (x_0)\not=0$, then we rotate the local
orthonormal frame about $x_0$ such that
\[
|v_1(x_0)|=|\nabla v(x_0)|\not=0\qquad \textrm{and}\qquad
v_i(x_0)=0,\quad i\geq2.
\]
From (\ref{max-prin}) we have at $x_0$,
\[
0=\frac12\nabla_{i} P=\sum_{j=1}^{n}v_jv_{ji}+Avv_i,
\]
\begin{equation}                \label{d1}
v_{11}=-Av \qquad \textrm{and}\qquad v_{1i}=0\quad i\geq 2,
\end{equation}
and
\begin{eqnarray}
&0 &\geq \frac12\Delta P(x_0) =\sum_{i, j=1}^{n}\left(v_{ji}v_{ji}+v_{j}v_{jii}+Av_{i}v_{i}+Avv_{ii}\right)\nonumber\\
&{} &=\sum_{i, j=1}^{n}\left(v_{ji}^2+v_j(v_{ii})_j +\textrm{R}_{ji}v_{j}v_{i}+Av_{ii}^2 +Av v_{ii}\right)\nonumber\\
&{} &=\sum_{i, j=1}^{n}v_{ji}^2+\nabla v\nabla(\Delta v) +\textrm{Ric}(\nabla v,\nabla v)+A|\nabla v|^2 +Av\Delta v\nonumber\\
&{} &\geq v_{11}^2+\nabla v\nabla(\Delta v) + (n-1)\kappa|\nabla v|^2+A|\nabla v|^2 +Av\Delta v\nonumber\\
&{} &=(-Av)^2-\lambda |\nabla v|^2+ (n-1)\kappa|\nabla v|^2+A|\nabla v|^2 -\lambda Av^2\nonumber\\
&{} &=[A-\lambda + (n-1)\kappa]|\nabla
v|^2+A(A-\lambda)v^2,\nonumber
\end{eqnarray}
where we have used (\ref{d1}) and (\ref{ricci-bound}). Therefore
at $x_0$,
\begin{equation}
0\geq [A-\lambda+(n-1)\kappa]|\nabla v|^2+A(A-\lambda)v^2.
\label{d2}
\end{equation}
That is,
\[
\epsilon |\nabla v(x_0)|^2 +
[-(n-1)\kappa+\epsilon][\lambda-(n-1)\kappa+\epsilon]v(x_0)^2\leq
0.
\]
Thus $\nabla v(x_0) = 0$. This contradicts the assumption $\nabla
v(x_0)\not= 0$.

Therefore we have $\nabla v(x_0)=0$, and
\[
P(x_0)=|\nabla v(x_0)|^2 +Av(x_0)^2=Av(x_0)^2\leq A.
\]
Now for all $x\in M$ we have
\[
|\nabla v(x)|^2+Av(x)^2=P(x)\leq P(x_0)\leq A
\]
and
\[
|\nabla v(x)|^2\leq A (1-v(x)^2).
 \]
Letting $\epsilon \rightarrow 0$ in the above inequality, the
estimate (\ref{basic3-d}) follows.
\end{proof}

We want to improve the above upper bound in (\ref{basic3-d})
further and proceed in the following way.

Define a function $Z$ on $[0,\sin^{-1}(1/b)]$ by
\[
Z(t)=\max_{x\in M,t=\sin^{-1} \left(v\left(x\right)/b\right)}
\frac{\left |\nabla v\right |^2}{b^2-v^2}/\lambda.
\]
The estimate in (\ref{basic3-d}) becomes
\begin{equation}                                \label{basic5-d}
Z(t)\leq 1+\beta\qquad \textrm{on } [0,\sin^{-1}(1/b)]
\end{equation}

Throughout this paper let
\[
\alpha=\frac12 (n-1)\kappa<0 \qquad \textrm{and} \qquad \delta =
\alpha/\lambda<0.
\]

We have the following conditions on the function $Z$.

\begin{theorem}                                                     \label{thm-barrier-d}
If the function $z:\ [0,\sin^{-1}(1/b)]\mapsto \mathbf{R}^1$
satisfies the following
\begin{enumerate}
 \item $z(t)\geq Z(t) \qquad t\in [0,\sin^{-1}(1/b)]$,
 \item there exists some $x_0\in M$
       such that at point $t_0=\sin^{-1} (v(x_0)/b)$ \linebreak
       $z(t_0)=Z(t_0)$,
 \item $z(t_0)\geq 1$,
 \item $z$ extends to a smooth even function, and
 \item $z'(t_0)\sin t_0\leq 0$,
\end{enumerate}
then we have the following
\begin{equation}                            \label{barrier-eq-d}
0\leq\frac12z''(t_0)\cos^2t_0 -z'(t_0)\cos t_0\sin t_0 -z(t_0)+
1-2\delta \cos^2t_0.
\end{equation}
\end{theorem}

\begin{proof}\quad
Define
\[ J(x)=\left\{ \frac{\left |\nabla v\right |^2}{b^2-v^2}
-\lambda z \right\}\cos^2t,
\]
where $t=\sin^{-1}(v(x)/b)$. Then
\[ J(x)\leq 0\quad\textrm{for } x\in M
\qquad \textrm{and} \qquad J(x_0)=0.
\]
So $J(x_0)$ is the maximum of $J$ on $M$. If $\nabla v(x_0)=0$,
then
\[ 0=J(x_0)=-\lambda z\cos^2 t.
\]
This contradicts the Condition 3 in the theorem. Therefore
\[ \nabla v(x_0)\not=0.
\]

We claim that
\begin{equation}\label{gra-j-eq-0}
\nabla J(x_0)=0.
\end{equation}
If $x_0\in M\backslash \partial M$, (\ref{gra-j-eq-0}) is
obviously true. Suppose that $x_0\in\partial M$. Take the same
local orthonormal frame $\{e_1, e_2,\cdots, e_{n}\}$ of $M$ about
$x_0$ as in the proof of Lemma \ref{pre-es-d-lemma}, where $e_n$
is the unit outward normal vector field near $x_0\in\partial M$,
$\{e_1, e_2,\cdots, e_{n-1}\}|_{\partial M}$ is a local frame of
$\partial M$ about $x_0$ and $\nabla_{e_n}e_i=0$ for $i\leq n-1$.
Since $v|_{\partial M}=0$, we have $v_i(x_0)=0$ for $i\leq n-1$.
$J(x_0)$ is the maximum implies that
\begin{equation}\label{j-i-eq-0}
J_i(x_0)= 0\qquad \textrm{for }i\leq n-1
\end{equation}
and
\begin{equation}\label{j-n-geq-0}
J_n(x_0)\geq 0.
\end{equation}
Using argument in proving (\ref{gra-v-leq-0}) and the
non-negativity of the mean curvature $m$ of $\partial M$ with
respect to the outward normal, we get
\[
\left(|\nabla v|^2\right)_n\Big|_{x_0}\leq 0.
\]
The Dirichlet condition $v(x_0)=0$ implies that $t(x_0)=0$ and
$z'(t(x_0))=z'(0)=0$, since by the Condition 4 in the theorem $z$
extends to a smooth even function. Therefore
\begin{equation}\label{j-n-leq-0}
J_n(x_0)=\frac{1}{b^2}\left(|\nabla v|^2\right)_n-\lambda\cos t[z'
\cos t -2z\sin t]t_n\Big|_{x_0}
 =\frac{1}{b^2}\left(|\nabla v|^2\right)_n\Big|_{x_0}
 \leq 0.
\end{equation}

Now (\ref{j-i-eq-0}), (\ref{j-n-geq-0}) and (\ref{j-n-leq-0})
imply (\ref{gra-j-eq-0}).

Thus (\ref{gra-j-eq-0}) holds, no matter $x_0\not\in\partial M$ or
$x_0\in\partial M$. By (\ref{gra-j-eq-0}) and the Maximum
Principle, we have
\begin{equation}                                                \label{es1}
\nabla J(x_0)=0\qquad \textrm{and}\qquad \Delta J(x_0)\leq 0.
\end{equation}
$J(x)$ can be rewritten as
\[  J(x)=\frac{1}{b^2}|\nabla v|^2-\lambda z\cos^2t.
\]
Thus (\ref{es1}) is equivalent to
\begin{equation}                                              \label{es2}
\frac{2}{b^2}\sum_{i}v_iv_{ij}\Big|_{x_0}=\lambda\cos t[z' \cos t
-2z\sin t]t_j\Big|_{x_0}
\end{equation}
and
\begin{eqnarray}                                              \label{es3}
0&\geq&\frac{2}{b^2}\sum_{i,j}v_{ij}^2+\frac{2}{b^2}\sum_{i,j}v_iv_{ijj}
 -\lambda (z''|\nabla t|^2+z'\Delta t)\cos^2t \\
 & &+4\lambda z'\cos t\sin t |\nabla t|^2 -
\lambda z\Delta\cos^2t\Big|_{x_0}.\nonumber
\end{eqnarray}
Rotate the frame so that $v_1(x_0)\not=0$ and $v_i(x_0)=0$ for
$i\geq 2$. Then (\ref{es2}) implies
\begin{equation}                                             \label{es4}
v_{11}\Big|_{x_0}=\frac{\lambda b}{2}(z'\cos t-2z\sin t)
\Big|_{x_0}\quad\text{and}\quad v_{1i} \Big|_{x_0}=0\ \text{for }
i\geq2.
\end{equation}
Now we have
\begin{eqnarray}
|\nabla v|^2
\Big|_{x_0}&=&\lambda b^2z\cos^2t\Big|_{x_0},\nonumber\\
 |\nabla t|^2
\Big|_{x_0}&=&\frac{|\nabla v|^2}{b^2-v^2}=\lambda z
\Big|_{x_0},\nonumber\\
\frac{\Delta v}{b}\Big|_{x_0} &=&\Delta \sin t =\cos t\Delta
t-\sin t |\nabla t|^2
\Big|_{x_0},\nonumber\\
\Delta t\Big|_{x_0}&=&\frac{1}{\cos t}(\sin t|\nabla
t|^2+\frac{\Delta v}{b})
\nonumber\\
 &=&\frac{1}{\cos t} [ \lambda z\sin t-\frac{\lambda}{b}v] \Big|_{x_0}, \quad\textrm{and}
\nonumber\\
\Delta\cos^2t\Big|_{x_0}&=&\Delta \left(1-\frac{v^2}{b^2}\right)
 =-\frac{2}{b^2}|\nabla v|^2-\frac{2}{b^2}v\Delta v
\nonumber\\
  &=&-2\lambda z\cos^2t+\frac{2}{b^2}\lambda v^2\Big|_{x_0}. \nonumber
\end{eqnarray}
Therefore,
\begin{eqnarray}
& {}&
\frac{2}{b^2}\sum_{i,j}v_{ij}^2\Big|_{x_0}\geq\frac{2}{b^2}v_{11}^2
\nonumber\\
& {}& =\frac{\lambda ^2}{2}(z')^2\cos^2t-2\lambda ^2zz'\cos t\sin
t
      +2\lambda ^2z^2\sin^2 t\Big|_{x_0}\nonumber,
\end{eqnarray}
\begin{eqnarray}
\frac{2}{b^2}\sum_{i,j}v_iv_{ijj}\Big|_{x_0}
&=&\frac{2}{b^2}\left(\nabla v\,\nabla
       (\Delta v)+\textrm{Ric}(\nabla v,\nabla v)\right)\nonumber\\
&\geq& \frac{2}{b^2}(\nabla v\,\nabla (\Delta v)+(n-1)\kappa|\nabla v|^2)\nonumber\\
 &=&-2\lambda^2z\cos^2t+4\alpha \lambda z\cos^2t\Big|_{x_0},\nonumber
\end{eqnarray}
\begin{eqnarray}
&{}&  -\lambda (z''|\nabla t|^2+ z'\Delta t)\cos^2t\Big|_{x_0}\nonumber\\
&{}&=-\lambda^2 zz''\cos^2t-
\lambda^2zz'\cos t\sin t\nonumber\\
&{}&{ }+\frac{1}{b}\lambda^2z'v\cos t \Big|_{x_0},\nonumber
\end{eqnarray}
and
\begin{eqnarray}
&{}&4\lambda z'\cos t\sin t|\nabla t|^2-\lambda z\Delta
\cos^2t\Big|_{x_0}
\nonumber\\
&{}&=4\lambda^2zz'\cos t\sin
t+2\lambda^2z^2\cos^2t-\frac{2}{b}\lambda^2zv\sin t
\Big|_{x_0}.\nonumber
\end{eqnarray}
Putting these results into (\ref{es3}) we get
\begin{eqnarray}                                                   \label{es5}
0&\geq&-\lambda^2zz''\cos^2t+ \frac{\lambda^2}{2}(z')^2\cos^2t
+\lambda^2z'\cos t\left(z\sin t +\sin t\right)
  \nonumber\\
 & & {}+2\lambda^2z^2-2\lambda^2z +4\alpha \lambda z\cos^2t
 \Big|_{x_0},
\end{eqnarray}
where we used (\ref{es4}). Now
\begin{equation}                                                    \label{es6}
z(t_0)>0,
\end{equation}
by Condition 3 in the theorem. Dividing two sides of (\ref{es5})
by $2\lambda^2z\Big|_{x_0}$, we have
\begin{eqnarray}
0&\geq&-\frac12z''(t_0)\cos^2t_0 +\frac12z'(t_0)\cos t_0\left(\sin t_0+\frac{\sin t_0}{z(t_0)}\right) +z(t_0) \nonumber\\
 & & {}  -1 +2\delta \cos^2t_0+\frac{1}{4z(t_0)}(z'(t_0))^2\cos^2t_0.\nonumber
\end{eqnarray}
Therefore,
\begin{eqnarray}
0&\geq&-\frac12z''(t_0)\cos^2t_0 + z'(t_0)\cos t_0\sin t_0+z(t_0) -1 +2\delta \cos^2t_0\nonumber\\
 & & {}+\frac{1}{4z(t_0)}(z'(t_0))^2\cos^2t_0+\frac12z'(t_0)\sin t_0\cos
 t_0[\frac{1}{z(t_0)}-1].                               \label{es6.1}
\end{eqnarray}
Conditions 2, 3 and 5 in the theorem imply that $z(t_0)=Z(t_0)\geq
1$ and $z'(t_0)\sin t_0\leq 0$. Thus the last two terms in
(\ref{es6.1}) are nonnegative and (\ref{barrier-eq-d}) follows.
\end{proof}

\section{Proof of Theorem \ref{main-thm}}\label{sec-proof}
\begin{proof}[Proof of Theorem  \ref{main-thm}]\quad
Let
\begin{equation}                                \label{z-def}
z(t)=1+\delta\xi(t),
\end{equation}
where $\xi$ is the functions defined by (\ref{xi-def}) in Lemma
\ref{xi-lemma}. We claim that
\begin{equation}                        \label{4.1-d}
Z(t)\leq z(t)\qquad \textrm{for }t\in  [0, \sin^{-1}(1/b)].
\end{equation}
Lemma \ref{xi-lemma} implies that for $t\in [0, \sin^{-1}(1/b)] $,
we have the following
\begin{eqnarray}
& &{}\frac{1}{2}z''\cos ^2t-z'\cos t\sin t-z
    =-1+ 2\delta\cos^2t,               \label{z-eq}\\
& &{}z'(t)\sin t\leq 0,\qquad (\textrm{since }\delta< 0)\label{z'-geq0}\\
& &{}z\textrm{ is a smooth even function, and}\\
& &{}z(t) \geq z(\frac{\pi}{2})=1. \label{z-min}
\end{eqnarray}

Let $P\in\mathbf{R}^1$ and $t_0\in [0,\sin^{-1}(1/b)]$ such that
\[ P=\max_{t\in [0,\sin^{-1}(1/b)]}\left(Z(t)-z(t)\right)=Z(t_0)-z(t_0).
\]
Thus
\begin{equation}\label{4.2}
Z(t)\leq z(t)+P\quad \textrm{for }t\in
[0,\sin^{-1}(1/b)]\qquad\textrm{and}\qquad Z(t_0)=z(t_0)+P.
\end{equation}
Suppose that $P>0$. Then $z+P$ satisfies the conditions in Theorem
\ref{thm-barrier-d} and therefore satisfies (\ref{barrier-eq-d}).
So we have
\begin{eqnarray}
&{}&z(t_0)+P=Z(t_0)\nonumber\\
&{}&\leq  \frac12(z+P)''(t_0)\cos^2 t_0-(z+P)'(t_0)\cos t_0
 \sin t_0+1-2\delta \cos^2 t_0\nonumber\\
&{}&=\frac12z''(t_0)\cos^2t_0-z'(t_0)\cos t_0\sin
t_0+1-2\delta \cos^2 t_0\nonumber\\
&{}&=z(t_0).\nonumber
\end{eqnarray}
This contradicts the assumption $P>0$. Thus $P\leq 0$ and
(\ref{4.1-d}) must hold. That along with the definition of the
function $Z$ means
\begin{equation}\label{4.3-d}
\sqrt{\lambda}\geq \frac{|\nabla t|}{\sqrt{z(t)}}.
\end{equation}

Take $q_1$ on $M$ such that $v(q_1)=1 =\sup_M v$ and and
$q_2\in\partial M$ such that distance $d(q_1, q_2) = \textrm{
distance }d(q_1, \partial M)$. Let $L$ be the minimum geodesic
segment between $q_1$ and $q_2$. We integrate both sides of
(\ref{4.3-d}) along $L$ and change variable and let $b\rightarrow
1$. Let $d$ be, as in Theorem \ref{main-thm}, the diameter of the
largest interior ball in $M$,$d=2r$ and $r=\max_{x\in M}
\textup{dist}(x,\partial M)$. Then
\begin{eqnarray}
\sqrt{\lambda}\,\frac{d}{2}&\geq& \int_{L}\,\frac{|\nabla
t(x)|}{\sqrt{z(t(x))}}\,dl= \int_0^{\frac{\pi}{2}}
\frac{1}{\sqrt{z(t)}}\,dt\nonumber\\
{}&\geq& \frac{\left(\int_0^{\pi/2}\
\,dt\right)^\frac32}{(\int_0^{\pi/2}\ z(t)\,dt)^{\frac12}} \geq
\left( \frac{(\frac{\pi}{2})^3}{\int_0^{\pi/2}\  z(t)\,dt}
\right)^{\frac12}\nonumber
\end{eqnarray}
Thus
\[
\lambda \geq \frac{\pi^3}{2d^2\int_0^{\pi/2} \ z(t)\,dt}.
\]
Now
\[
\int_0^{\frac{\pi}{2}}\ z(t)\,dt=\int_0^{\frac{\pi}{2}}\ [1+
\delta \xi(t)]\,dt=\frac{\pi}{2}(1-\delta),
\]
by (\ref{xi-int}) in Lemma \ref{xi-lemma}. Therefore we have
\[
\lambda \geq \frac{\pi^2}{(1-\delta)d^2}\quad\textrm{and}\quad
\lambda \geq  \frac12(n-1)\kappa+ \frac{\pi^2}{d^2}.
\]
\end{proof}

We now present a lemma that is used in the proof of Theorem
\ref{main-thm}.

\begin{lemma}                           \label{xi-lemma}
Let
\begin{equation}                        \label{xi-def}
\xi(t)=\frac{\cos^2t+2t\sin t\cos t +t^2-\frac{\pi^2}{4}}{\cos^2t}
\qquad \textrm{on}\quad [-\frac{\pi}{2},\frac{\pi}{2}\,].
\end{equation}
Then the function $\xi$  satisfies the following
\begin{eqnarray}
& &{}\frac{1}{2}\xi''\cos ^2t-\xi'\cos t\sin t-\xi
    =2\cos^2t\quad \textrm{in }(-\frac{\pi}{2},\frac{\pi}{2}\,),          \label{xi-eq}\\
& &{}\xi'\cos t -2\xi\sin t =4t\cos t                      \label{xi-eq2}\\
& &{}\int_0^{\frac{\pi}{2}}\xi(t)\, dt= -\frac{\pi}{2}          \label{xi-int}\\
& &{}1-\frac{\pi^2}{4}=\xi(0)\leq\xi(t)\leq\xi(\pm
\frac{\pi}{2})=0\quad
\textrm{on }[-\frac{\pi}{2},\frac{\pi}{2}\,],                                \nonumber\\
& &{} \xi' \textrm{ is increasing on }
[-\frac{\pi}{2},\frac{\pi}{2}\,] \textrm{ and }
\xi'(\pm \frac{\pi}{2}) =\pm \frac{2\pi}{3},                     \nonumber\\
& &{}\xi'(t)< 0 \textrm{ on }(-\frac{\pi}{2},0)\textrm{ and \ }
\xi'(t)>0 \textrm{ on }(0,\frac{\pi}{2}\,), \nonumber \\
& &{}\xi''(\pm\frac{\pi}{2})=2, \ \xi''(0)=2(3-\frac{\pi^2}{4})
\textrm{ and \ } \xi''(t)> 0 \textrm{ on }
[-\frac{\pi}{2},\frac{\pi}{2}\,], \nonumber\\
& &{}(\frac{\xi'(t)}{t})'>0 \textrm{ on } (0,\pi/2\,)\textrm{ and
\ } 2(3-\frac{\pi^2}{4})\leq \frac{\xi'(t)}{t}\leq \frac43
\textrm{ on } [-\frac{\pi}{2},\frac{\pi}{2}\,],\nonumber \\
& &{}\xi'''(\frac{\pi}{2})=\frac{8\pi}{15}, \xi'''(t)< 0 \textrm{
on }(-\frac{\pi}{2},0) \textrm{ and \ }  \xi'''(t)>0 \textrm{ on
}(0,\frac{\pi}{2}\,). \nonumber
\end{eqnarray}
\end{lemma}
\begin{proof}\quad
For convenience, let $q(t)= \xi'(t)$, i.e.,
\begin{equation}                                       \label{q-def}
q(t) = \xi'(t) = \frac{2(2t\cos t +t^2\sin t +\cos^2 t \sin t
-\frac{\pi^2}{4}\sin t)}{\cos^3 t}.
\end{equation}
Equation (\ref{xi-eq}) and the values $\xi(\pm \frac{\pi}{2})=0$,
$\xi(0)=1-\frac{\pi^2}{4}$ and $\xi'(\pm \frac{\pi}{2}) =\pm
\frac{2\pi}{3}$ can be verified directly from (\ref{xi-def}) and
(\ref{q-def}) .  The values of $\xi''$ at $0$ and $\pm
\frac{\pi}{2}$ can be computed via (\ref{xi-eq}). By
(\ref{xi-eq2}), $(\xi(t)\cos^2 t)' =4t\cos^2 t$. Therefore
\newline $\xi(t)\cos^2 t=\int_{\frac{\pi}{2}}^t \ 4s\cos^2 s\,ds$,
and
\[
\int_{-\frac{\pi}{2}}^{\frac{\pi}{2}}\
\xi(t)\,dt=2\int_0^{\frac{\pi}{2}}\
\xi(t)\,dt=-8\int_0^{\frac{\pi}{2}}\left( \frac{1}{\cos^2(t)}
\int_t^{\frac{\pi}{2}}\ s\cos^2s\,ds\right)\,dt
\]
\[
=-8\int_0^{\frac{\pi}{2}}\left(\int_0^s\
\frac{1}{\cos^2(t)}\,dt\right)\ s\cos^2s\,ds
=-8\int_0^{\frac{\pi}{2}}\ s\cos s\sin s\,ds=-\pi.
\]
It is easy to see that $q$ and $q'$ satisfy the following
equations
\begin{equation}                                         \label{q-eq}
\frac12 q''\cos t -2q'\sin t -2q\cos t = -4 \sin t,
\end{equation}
and
\begin{equation}                                       \label{q'-eq}
\frac{\cos^2 t}{2(1+\cos^2 t)}(q')''-\frac{2\cos t\sin t}{1+\cos^2
t}(q')'-2(q')=-\frac{4}{1+\cos^2 t}.
\end{equation}
The last equation implies $q'=\xi''$ cannot achieve its
non-positive local minimum at a point in $(-\frac{\pi}{2},
\frac{\pi}{2})$. On the other hand, $\xi''(\pm\frac{\pi}{2})=2$,
by equation (\ref{xi-eq}), $\xi(\pm \frac{\pi}{2})=0$ and
$\xi'(\pm \frac{\pi}{2})=\pm \frac{2\pi}{3}$. Therefore
$\xi''(t)>0$ on $[-\frac{\pi}{2},\frac{\pi}{2}]$ and $\xi'$ is
increasing. Since $\xi'(t)=0$, we have $\xi'(t)< 0$ on
$(-\frac{\pi}{2},0)$ and $\xi'(t)>0$ on $(0,\frac{\pi}{2})$.
Similarly, from the equation
\begin{eqnarray}                                           \label{q''-eq}
&\frac{\cos^2 t}{2(1+\cos^2 t)}(q'')'' -\frac{\cos t\sin t
(3+2\cos^2 t)}{(1+\cos^2 t)^2}(q'')' -\frac{2(5\cos^2 t+\cos^4 t)}{(1+\cos^2 t)^2}(q'') \nonumber\\
&=-\frac{8\cos t\sin t}{(1+\cos^2 t)^2}
\end{eqnarray}
we get the results in the last line of the lemma.

Set $h(t)=\xi''(t)t-\xi'(t)$. Then $h(0)=0$  and $h'(t)=
\xi'''(t)t>0$ in $(0,\frac{\pi}{2})$. Therefore
$(\frac{\xi'(t)}{t})'=\frac{h(t)}{t^2}>0$ in $(0,\frac{\pi}{2})$.
Note that $\frac{\xi'(-t)}{-t}= \frac{\xi'(t)}{t}$,
$\frac{\xi'(t)}{t}|_{t=0}=\xi''(0)=2(3-\frac{\pi^2}{4})$ and
$\frac{\xi'(t)}{t}|_{t=\pi/2}=\frac43$. This completes the proof
of the lemma.
\end{proof}

%
%

Department of mathematics, Utah Valley State College, Orem, Utah 84058

\textit {E-mail address}: \texttt{lingju@uvsc.edu}
\end{document}